\documentclass[12pt]{amsart}
\usepackage{amsfonts,amsthm,latexsym,amsmath,amssymb, 
amscd,amsmath, mathrsfs, url, cite, filecontents, mathtools}
\newcounter{count}

\numberwithin{count}{section}

\newtheorem{Corollary}[count]{Corollary}
\newtheorem{Theorem}[count]{Theorem}

\begin{document}

\author[K.~Bielenova ]{Kateryna Bielenova }

\address{Department of Mathematics \& Computer Sciences, 
V. N. Karazin Kharkiv National University,
4 Svobody Sq., Kharkiv, 61022, Ukraine}
\email{belenooova@gmail.com}

\author[H.~Nazarenko]{Hryhorii Nazarenko}

\address{Department of Mathematics \& Computer Sciences, 
V. N. Karazin Kharkiv National University,
4 Svobody Sq., Kharkiv, 61022, Ukraine}
\email{548gggg@gmail.com}

\author[A.~Vishnyakova]{Anna Vishnyakova}
\address{Department of Mathematics \& Computer Sciences, 
V. N. Karazin Kharkiv National University,
4 Svobody Sq., Kharkiv, 61022, Ukraine}
\email{anna.vishnyakova@karazin.ua}

\title[Complex polynomials with simple zeros]
{A sufficient condition for a complex polynomial to have only simple 
zeros and an analog of Hutchinson's theorem for real polynomials}

\begin{abstract}

We find the constant $b_{\infty}$ ($b_{\infty} \approx 4.81058280$)
such that if a complex polynomial or entire function $f(z) = \sum_{k=0}^
\omega  a_k z^k, $ $\omega \in \{2, 3, 4, \ldots  \} \cup \{\infty\},$ with
nonzero coefficients satisfy the conditions $\left|\frac{a_k^2}{a_{k-1} 
a_{k+1}}\right| >b_{\infty} $ for all $k =1, 2, \ldots, \omega-1,$ then 
all the zeros of $f$ are simple. We show that the constant $b_{\infty}$ in 
the statement above is the  smallest possible. We also obtain an analog of  
Hutchinson's theorem for polynomials or entire functions with real nonzero 
coefficients.

\end{abstract}



\maketitle

\section{Introduction}

In this short note, we obtain a  simple sufficient condition for a complex 
polynomial to have only simple zeros in terms of its coefficients. To 
formulate our results, we define the second quotients of coefficients for 
a polynomial.

Let us consider a complex polynomial (or entire function) $f(z)=
\sum_{k=0}^{\omega }a_{k}z^k$, where $a_{k}\in \mathbb{C} \setminus \{0\}$ 
and $\omega \in \{2, 3, 4, \ldots  \} \cup \{\infty\}.$ 
We define  the second quotients of the Taylor coefficients of $f$ 
by the formula
\begin{equation} 
\label{f1}
q_{n}(f) =\frac{a_{n-1}^2}{a_{n-2}a_{n}} , \   n\geq 2.
\end{equation}
It is easy to check that
\begin{equation} 
\label{f2}
a_{n}=\frac{a_1^n}{a_0^{n-1} q_{2}^{n-1}q_{3}^{n-2}
\cdot \ldots \cdot q_{n-1}^2q_{n}}, \    n \geq 2. 
\end{equation}
One can see that the second quotients of Taylor coefficients are 
independent parameters that define a function up to multiplication 
by a constant and changing $z$ to $\lambda z.$

In 1926, J.I.~Hutchinson found quite a simple sufficient condition 
for an entire function with positive coefficients to have only real simple zeros.

{\bf  Theorem A} (J.I.~Hutchinson, \cite{hut}). { \it  Let 
$f(z)=\sum_{k=0}^{\infty}a_{k}z^k, a_{k}>0$ for all k, be an entire function. Then the 
inequalities $q_{n}(f) \geq 4$ for all $n \geq 2$ hold
if and only if the following two conditions are fulfilled:

(i) The zeros of $f$ are all real, simple and negative, and

(ii) the zeros of any polynomial $\sum_{k=m}^{n}a_{k}z^k,  m < n, $ 
formed by taking any number of consecutive terms of $f,$ are all 
real and non-positive.}

For some extensions of Hutchinson's results see, for example, \cite{cc1}, where, 
in particular, the following theorem is proved.

{\bf  Theorem B} (T.~Craven and G.~Csordas, \cite{cc1}).  { \it   Let $N\in \mathbb{N}$ and 
$(\gamma_k)_{k=0}^N,$ $\gamma_0 =1,$  be a sequence of positive real numbers. Suppose 
that the inequalities  $ \frac{\gamma_n^2}{\gamma_{n-1} \gamma_{n+1}} \geq \alpha^2$  
hold for all $n=1,2, \ldots, N-1$, where $\alpha = \max(2, \frac{\sqrt{2}}{2}
\left(1+ \sqrt{1+\gamma_1}  \right).$  Then the polynomial  $Q(x) =\sum_{n=0}^N \gamma_n \cdot 
\frac{x (x-1) \cdot \ldots \cdot (x-n+1)}{n!}$   has only real, simple, negative zeros.}

There are a number of works which deal with the statements of the following kind: there 
exists a constant $d>1$ such that if a real polynomial $P$ satisfies the condition $q_k(P) 
> d$ for all $k$, then we can state something about the location of zeros  $P.$ For 
example, in \cite{hen} the author proved that  if  for some constant $d>0$ a real 
polynomial $P$ satisfies the condition $q_k(P) > d$ for all $k,$ then all the
zeros of $P$ lie in a special sector depending on $d.$   In \cite{KatVi} the 
smallest possible constant $d>0$ was found
such that  if a real polynomial $P$ satisfies the condition $q_k(P) > d$ for all $k,$
then $P$ is stable (all the zeros of $P$ lie in the left half-plane). In this paper, we 
study analogous questions for complex polynomials and entire functions.

Hutchinson's theorem inspired our investigations. The goal of this work is to
find  sufficient conditions for complex polynomials or entire functions with
non-zero coefficients to have only simple zeros. More precisely, we answer the 
following question:  what is the smallest possible constant $c > 0$ such that for 
every complex polynomial $P$  with nonzero coefficients 
if the inequalities $q_{n}(P) > c$  hold for all $n \geq 2,$ then all 
the zeros of $P$ are simple.

For $x>1$ let us consider the function $\phi(x)=1-2\sum_{k=1}^{\infty}
x^{\frac{- k^2}{2}}.$  
We observe  that $\phi$  is an increasing function on  $(0;\infty), $  
$\lim_{x \to 1+0} \phi(x) = 
-\infty$ and $\lim_{x \to +\infty} \phi(x)=1.$  So, the equation
\begin{equation} 
\label{f3}
1-2\sum_{k=1}^{\infty}x^{\frac{- k^2}{2}}=0
\end{equation}
has a unique positive root, which we denote by $b_{\infty}$. 
One can check that $b_{\infty} \approx 4.81058280$.

For $n \in \mathbb{N}$  we also define $b_{2n}$ as the unique positive 
root of the equation
\begin{equation} 
\label{f4}
1-2\sum_{k=1}^{n} x^\frac{-k^2}{2}=0.
\end{equation}
One can see that $\left (b_{2n} \right)_{n=1}^\infty$ is the
increasing sequence, $\lim_{n\to\infty}  b_{2n} = b_{\infty},$ and
\begin{equation} 
\label{f5}
b_2 = 4,  b_{4} \approx 4{.}79753651.
\end{equation}

The constant  $b_{\infty}$ firstly appeared in the paper \cite{KarVi}
where some analogs of Hutchinson's result were obtained for sign-independently 
hyperbolic polynomials. 

{\bf  Theorem C} (I. Karpenko and A. Vishnyakova, \cite{KarVi}). { \it  Let $f(z)= 
\sum_{k=0}^{\infty}a_{k}z^k$ be an entire function with positive coefficients. 
Suppose that $q_{k}(f) \geq b_{\infty}$ for all $k\geq2$. Then for every $n \in \mathbb{N},  $ 
the $n$-th section $S_{n}(z) :=\sum_{k=0}^{n}a_{k}z^k$ is a sign-independently  sign-hyperbolic polynomial, meaning that it remains real-rooted after arbitrary sign changes 
of its coefficients. }

Our first result is the following theorem.

\begin{Theorem}
\label{th:mthm1}
(i) Let $n \in \mathbb{N} $ be a given integer, and $P_{2n}(z) = 
\sum_{k=0}^{2n} a_k z^k $, $a_k 
\in \mathbb{C} \setminus \{0\}$ for all $k,$  be a polynomial. 
Suppose that the inequalities  
$|q_k(P_{2n})| > b_{2n}$ hold for all $k = 2, 3, \ldots , 2n.$ 
Then all the zeros of $P_{2n}$ are simple. Moreover, the moduli of 
all zeros of $P_{2n}$  are pairwise different.

(ii) Let $n \in \mathbb{N} $ be a given integer, and $P_{2n+1}(z) = 
\sum_{k=0}^{2n+1} a_k z^k $, $a_k 
\in \mathbb{C} \setminus \{0\}$ for all $k,$  be a polynomial. 
Suppose that the inequalities  
$|q_k(P_{2n+1})| \geq b_{2n+2}$ hold for all $k = 2, 3, \ldots , 2n+1.$ 
Then all the zeros of $P_{2n+1}$ are simple. Moreover, the moduli of 
all zeros of $P_{2n+1}$  are pairwise different.

(iii)  Let $f(z)=\sum_{k=0}^\infty a_k z^k $, $a_k \in \mathbb{C} 
\setminus \{0\}$ for all $k,$  be an 
entire function.  Suppose that the inequalities  
$|q_k(f)| \geq b_{\infty}$ hold for all $k \geq 2.$ 
Then all the zeros of $f$ are simple. Moreover, the moduli of 
all zeros of $f$  are pairwise different.
\end{Theorem}

Since the sequence $(b_{2n})_{n=1}^\infty$ is monotonic
and tends to $b_\infty,$ we get the following corollary.

\begin{Corollary}
\label{th:cor1}
Let $n \geq 2 $ be a given integer, and $P$  be a 
complex polynomial with nonzero coefficients of degree $n.$
If the inequalities  $|q_k(P)| \geq b_{\infty}$ hold for all 
$k = 2, 3, \ldots , n,$  then all the zeros of $P$ are simple.
Moreover, the moduli of all zeros of $P$  are pairwise different.
\end{Corollary}

The following statement shows the sharpness of Theorem \ref{th:mthm1} 
for entire functions and polynomials of even degrees, and asymptotical 
sharpness of Theorem \ref{th:mthm1} for polynomials of odd degrees.

\begin{Theorem}
\label{th:mthm2}
(i) For every $n \in \mathbb{N} $ there exists a complex polynomial 
$P_{2n}(z) = \sum_{k=0}^{2n} a_k z^k $, $a_k 
\in \mathbb{C} \setminus \{0\}$ for all $k,$  such that   the  equalities  
$|q_k(P_{2n})| = b_{2n}$ hold for all $k = 2, 3, \ldots , 2n,$  and  
 $P_{2n}$ has a multiple  root.
 
(ii) For every $n \in \mathbb{N} $ and every $\varepsilon >0 $
there exists a complex polynomial   $P_{2n+1, \varepsilon }$ 
with nonzero coefficients, $\deg P_{2n+1, \varepsilon } = 2n+1,$
such that  $|q_k(P_{2n+1, \varepsilon })| > b_{2n} -
\varepsilon $  for all $k = 2, 3, \ldots , 2n+1,$
and $P_{2n+1, \varepsilon }$ has a multiple  root. 
 
(iii)  For every $\varepsilon >0$ there exists an entire function 
$f_\varepsilon(z)=\sum_{k=0}^\infty a_k(\varepsilon) z^k $, 
$a_k(\varepsilon) \in \mathbb{C} \setminus \{0\}$ for all $k,$  
such that the inequalities  $|q_k(f_\varepsilon)| > b_{\infty} 
- \varepsilon$ hold for all $k \geq 2,$ and  $f$ has a multiple  root.

\end{Theorem}

In the following statement we find the sharp constant
for polynomials of the third degree.

\begin{Theorem}
\label{th:mthm3}
(i) Let $P_3$  be  a complex polynomial 
with nonzero coefficients, $\deg P_{3 } = 3,$
and suppose that  $|q_k(P_{3 })| > \sqrt{9+6\sqrt{3}}  $  
for  $k = 2, 3.$ Then all the zeros of $P_{3}$ are simple.
Moreover, the moduli of all zeros of $P_{3}$  are pairwise different.
Note that $ \sqrt{9+6\sqrt{3}}\approx  4{.}4036695,$ 
and  $4= q_2 < \sqrt{9+6\sqrt{3}} < q_4 
\approx 4{.}79753651.$
 
(ii) There exists a complex polynomial   $Q_{3}$ 
with nonzero coefficients, $\deg Q_{3 } = 3,$
such that  $|q_k(Q_{3})| = \sqrt{9+6\sqrt{3}}$  for 
$k = 2, 3,$ and $Q_{3}$ has a multiple  root. 
 \end{Theorem}
 
Using Theorem \ref{th:mthm1} and Theorem \ref{th:mthm3} (i) 
we obtain the following analog of Hutchinson's theorem for real 
polynomials.
 
 \begin{Theorem}
\label{th:mthm4}
(i)  Let $n \in \mathbb{N} $ be a given integer, and $P_{2n}(z) = 
\sum_{k=0}^{2n} a_k z^k $, $a_k 
\in \mathbb{R} \setminus \{0\}$ for all $k,$  be a real polynomial. 
Suppose that the inequalities  
$|q_k(P_{2n})| \geq b_{2n}$ hold for all $k = 2, 3, \ldots , 2n.$ 
Then all the zeros of $P_{2n}$ are real.

(ii)  For every $n \in \mathbb{N} $ and every $\varepsilon >0 $
there exists a real polynomial   $P_{2n, \varepsilon }$ 
with nonzero coefficients, $\deg P_{2n, \varepsilon } = 2n,$
such that  $|q_k(P_{2n, \varepsilon })| > b_{2n} -
\varepsilon $  for all $k = 2, 3, \ldots , 2n,$
and $P_{2n, \varepsilon }$  has nonreal roots.  

(iii)  Let $n \in \mathbb{N} $ be a given integer, and $P_{2n+1}(z) = 
\sum_{k=0}^{2n+1} a_k z^k $, $a_k 
\in \mathbb{R} \setminus \{0\}$ for all $k,$  be a real polynomial. 
Suppose that the inequalities  
$|q_k(P_{2n+1})| \geq b_{2n+2}$ hold for all $k = 2, 3, \ldots , 2n+1.$ 
Then all the zeros of $P_{2n+1}$ are real.

(iv)  For every $n \in \mathbb{N} $ and every $\varepsilon >0 $
there exists a real polynomial   $P_{2n+1, \varepsilon }$ 
with nonzero coefficients, $\deg P_{2n+1, \varepsilon } = 2n+1,$
such that  $|q_k(P_{2n+1, \varepsilon })| > b_{2n} -
\varepsilon $  for all $k = 2, 3, \ldots , 2n+1,$
and $P_{2n+1, \varepsilon }$ has nonreal roots.

(v) Let $P_3$  be  a real polynomial 
with nonzero coefficients, $\deg P_{3 } = 3,$
and suppose that  $|q_k(P_{3 })| \geq \sqrt{9+6\sqrt{3}}  $  
for  $k = 2, 3.$ Then all the zeros of $P_{3}$ are real.

(vi) For every $\varepsilon >0 $
there exists a real polynomial   $P_{3, \varepsilon }$ 
with nonzero coefficients, $\deg P_{3, \varepsilon } = 3,$
such that  $|q_k(P_{3, \varepsilon })| > \sqrt{9+6\sqrt{3}} -
\varepsilon $  for  $k = 2, 3,$
and $P_{3, \varepsilon }$  has nonreal roots. 

 \end{Theorem}

We see that, unlike Hutchinson's result, the sharp constant
for the realrootedness of a real polynomial depends on the
degree of the polynomial. 

\section{Proof of  Theorem  \ref{th:mthm1}}

At first we consider the case of polynomials of even degrees. Let $n \in \mathbb{N}, $ 
$P_{2n}(z)=\sum_{k=0}^{2n}a_{k}z^k$, where $a_{k}\in \mathbb{C} 
\setminus \{0\},$ and suppose that the inequalities  
$|q_k(P_{2n})| > b_{2n}$ hold for all $k = 2, 3, \ldots, 2n.$ 

Without loss of generality, we can assume that $a_0=a_1=1,$ since we can 
consider the function $Q_{2n}(z) =a_0^{-1} P_{2n} (a_0 a_1^{-1} z) $  instead of 
$P_{2n},$ due to the fact that such a rescaling of $P_{2n}$ preserves its property of 
having all simple zeros  and preserves the second quotients:  $q_k(Q_{2n}) =q_k(P_{2n})$ 
for all $k.$ During the proof we use the notation  $q_k$ instead of 
 $q_k(P_{2n}).$ So, we can write 
$$P_{2n}(z)=1 + z+\frac{z^2}{q_{2}}+\frac{z^3}{q_{2}^2q_{3}}
+\ldots
+ \frac{z^{2n}}{q_{2}^{2n-1}q_{3}^{2n-2}\cdot \ldots \cdot  q_{2n-1}^2q_{2n}}.$$ 

We choose an arbitrary $ k= 1, 2, 3, \ldots, 2n -1$ and fix this $ k.$  
Denote by 
\begin{equation} 
\label{f00}
R_1 := \sqrt{|q_{2}|}, \   R_{k}:=|q_{2}q_{3}\dots q_{k}\sqrt{q_{k+1}}|, \  
 k=  2, 3, \ldots, 2n -1. 
\end{equation}  
We have 

$$P_{2n}(z)=\sum_{j=0}^{k-1}\frac{z^j}{q_{2}^{j-1}q_{3}^{j-2}\cdot \ldots \cdot  q_{j}}+
\frac{z^k}{q_{2}^{k-1}q_{3}^{k-2} \cdot \ldots \cdot q_{k}}+$$
$$\sum_{j=k+1}^{2n}
\frac{z^j}{q_{2}^{j-1}q_{3}^{j-2}\cdot \ldots \cdot  q_{j}} 
 =:  S_{1, k}(z)+S_{2, k}(z)+S_{3, k} (z). $$
 
We want to prove  the inequality 
\begin{equation} 
\label{f60}\min_{|z| = R_k}|S_{2, k} (z)| > 
\max_{|z| = R_k}( |S_{1, k}(z)|+|S_{3, k} (z)|).
\end{equation}  
We obtain
for every $z, |z| = R_k,$

\begin{equation} 
\label{f6}
|S_{2, k} (z)|=|q_{2}q_{3}^2 \cdot \ldots \cdot  q_{k}^{k-1}q_{k+1}^{k/2}|.
\end{equation}

Now we estimate from above $|S_{1, k}(z)|$ for $|z|=R_{k}.$
We have

$$  |S_{1, k}(z)| \leq
\sum_{j=0}^{k-1}\frac{|z^j|}{|q_{2}^{j-1}q_{3}^{j-2}  \cdot \ldots \cdot q_{j}|}=
\sum_{j=0}^{k-1}|
\frac{q_{2}^{j}q_{3}^j \cdot \ldots \cdot q_{k}^j q_{k+1}^{j/2}}{q_{2}^{j-1}
q_{3}^{j-2} \cdot \ldots \cdot q_{j}}|=$$
$$=|q_{2} q_{3}^2  \cdot \ldots \cdot q_{k-1}^{k-2}q_{k}^{k-1}q_{k+1}^{(k-1)/2}|+|q_{2}  
q_{3}^2 \cdot \ldots \cdot  q_{k-2}^{k-3}q_{k-1}^{k-2}q_{k}^{k-2}q_{k+1}^{(k-2)/2}
|$$
$$ + |q_{2}  
q_{3}^2 \cdot \ldots \cdot  q_{k-2}^{k-3}q_{k-1}^{k-3}q_{k}^{k-3}q_{k+1}^{(k-3)/2}
| + \ldots  + |q_{2}q_{3} \cdot \ldots \cdot q_{k}\sqrt{q_{k+1}}| + 1    $$ 
(we rewrite the sum from the end to the beginning). Thus, we get
\begin{eqnarray} 
\nonumber
&  |S_{1, k} (z)|\leq|q_{2} q_{3}^2  \cdot \ldots \cdot 
q_{k-1}^{k-2}q_{k}^{k-1}q_{k+1}^{k/2}|
\cdot \left(\left|\frac{1}{q_{k+1}^{1/2}}\right| + \left|\frac{1}{q_{k}
q_{k+1}}\right|+ \right.
\\
\nonumber  & 
\left|\frac{1}{q_{k-1}q_{k}^2q_{k+1}^{3/2}}\right| + \ldots + 
\left|\frac{1}{q_{k-j+1} q_{k-j+2}^2 \cdot \ldots \cdot 
q_{k-1 }^{j-1} q_k^j q_{k+1}^{(j+1)/2}}\right| 
\\
\nonumber  & \left. + \ldots + \left|\frac{1}{q_3 q_4^2 \cdot \ldots \cdot 
q_k^{k-2}q_{k+1}^{(k-1)/2}}\right|  
+ \left|\frac{1}{q_2 q_3^2 q_4^3 \cdot \ldots \cdot 
q_k^{k-1}q_{k+1}^{k/2}}\right| \right). 
\end{eqnarray}

Using our assumption $|q_k(P_{2n})| > b_{2n}$ for $k = 2, 3, \ldots, 2n,$
we obtain
\begin{equation} 
\label{f8}
|S_{1, k} (z)| < \left|q_{2} q_{3}^2  \cdot \ldots \cdot q_{k-1}^{k-2}
q_{k}^{k-1}q_{k+1}^{k/2}\right| \cdot \sum_{j=1}^{k}b_{2n}^\frac{-j^2}{2}.
\end{equation}

Now we estimate $|S_{3,k}(z)|$ from above for $|z|=R_{k}.$
We have
$$ |S_{3,k}(z)|
\leq  \sum_{j=k+1}^{2n}\left |\frac{q_{2}^{j}q_{3}^j
\cdot \ldots \cdot  q_{k}^j q_{k+1}^{j/2}}{q_{2}^{j-1}q_{3}^{j-2}
\cdot \ldots \cdot q_{j}}\right|=$$
$$=|q_{2} q_3^2 \cdot \ldots \cdot q_{k}^{k-1}q_{k+1}^{(k-1)/2}|+\left|\frac{q_{2} q_3^2
 \cdot \ldots \cdot q_{k}^{k-1}
q_{k+1}^{(k-2)/2}}{q_{k+2}}\right|+$$
$$+ \ldots + \left|\frac{q_{2} q_3^2
 \cdot \ldots \cdot q_{k}^{k-1}
q_{k+1}^{(2k- j)/2}}{q_{k+2}^{j -k-1}\cdot \ldots \cdot q_{j-1}^2 q_{j}}\right|
+ \ldots + \left|\frac{q_{2} q_3^2
 \cdot \ldots \cdot q_{k}^{k-1}
q_{k+1}^{(2k-2n)/2}}{q_{k+2}^{2n-k-1}\cdot \ldots \cdot q_{2n-1}^2 q_{2n}}\right|  $$
$$= \left|q_{2} q_{3}^2  \cdot \ldots \cdot q_{k-1}^{k-2}
q_{k}^{k-1}q_{k+1}^{k/2}\right| \left( \left|\frac{1}{q_{k+1}^{1/2}}\right|+
\left|\frac{1}{q_{k+1}q_{k+2}}\right|+\ldots +\right.$$
$$ \left. \left|\frac{1}{q_{k+1}^{(j-k)/2}q_{k+2}^{j -k-1}
\cdot \ldots \cdot q_{j-1}^2 q_{j}}\right|
+ \ldots + \left|\frac{1}{q_{k+1}^{(2n-k)/2}q_{k+2}^{2n-k-1}
\cdot \ldots \cdot q_{2n-1}^2 q_{2n}}\right|     \right).    $$

Using our assumption $|q_k(P_{2n})| > b_{2n}$ for $k = 2, 3, \ldots, 2n,$
we obtain
\begin{equation} 
\label{f9}
|S_{3, k} (z)| < \left|q_{2} q_{3}^2  \cdot \ldots \cdot q_{k-1}^{k-2}
q_{k}^{k-1}q_{k+1}^{k/2}\right| \cdot \sum_{j=1}^{2n-k}b_{2n}^\frac{-j^2}{2}.
\end{equation}
Thus, by virtue of  (\ref{f6}), (\ref{f8}) and (\ref{f9}),  the desired 
inequality (\ref{f60}) follows from
\begin{equation} 
\label{f10}
1-\sum_{j=1}^{k}b_{2n}^\frac{-j^2}{2} - \sum_{j=1}^{2n-k}
b_{2n}^\frac{-j^2}{2}\geq 0.
\end{equation}
Since the summands in both sums in the above inequality are strictly
decreasing in $j$ we have
\begin{equation} 
\label{f11}
1-\sum_{j=1}^{k}b_{2n}^\frac{-j^2}{2} - \sum_{j=1}^{2n-k}
b_{2n}^\frac{-j^2}{2}\geq 1-2\sum_{j=1}^{n}
b_{2n}^\frac{-j^2}{2}=0
\end{equation}
by the definition (\ref{f4}) of the constant $b_{2n}.$
We have proved that for every $ k= 1, 2, 3, \ldots, 2n -1$
the inequality (\ref{f60}) is valid.  Thus, by  Rouch\'{e}'s 
theorem, we obtain that for every $ k= 1, 2, 3, \ldots, 2n -1$
the polynomial $P_{2n}$ has exactly $k$ zeros in the circle
$\{z : |z| < R_k  \}.$  Whence, the polynomial $P_{2n}$
has one zero in the circle $\{z : |z| < R_1  \},$ one
zero in the annulus $\{z : R_1 \leq |z| < R_2  \},$
one zero in the annulus $\{z : R_2 \leq |z| < R_3  \},$
and so on, one zero in the annulus $\{z : R_{2n-2} \leq |z| < R_{2n-1}  \}$
and one zero in the set $\{z :  |z| \geq R_{2n-1}  \}.$
We have proved that all the zeros of $P_{2n}$ are simple.
Moreover, the moduli of all zeros of $P_{2n}$  are pairwise different.

Now we consider the case of polynomials of odd degrees.
Let $n\in \mathbb{N}$  and
$$P_{2n+1}(z)=1 + z+\frac{z^2}{q_{2}}+\frac{z^3}{q_{2}^2q_{3}}
+\ldots
+ \frac{z^{2n+1}}{q_{2}^{2n}q_{3}^{2n-1}\cdot \ldots \cdot  q_{2n}^2q_{2n+1}}$$ 
be a complex polynomial with nonzero coefficients. Suppose that the inequalities  
$|q_k(P_{2n+1})| \geq b_{2n+2}$ hold for all $k = 2, 3, \ldots , 2n+1.$  We choose an arbitrary 
$ k= 1, 2, 3, \ldots, 2n $ and fix this $ k.$  We use the same notation $R_k$ as in
(\ref{f00}).
We have 
$$P_{2n+1}(z)=\sum_{j=0}^{k-1}\frac{z^j}{q_{2}^{j-1}q_{3}^{j-2}\cdot \ldots \cdot  q_{j}}+
\frac{z^k}{q_{2}^{k-1}q_{3}^{k-2} \cdot \ldots \cdot q_{k}}+$$
$$\sum_{j=k+1}^{2n+1}
\frac{z^j}{q_{2}^{j-1}q_{3}^{j-2}\cdot \ldots \cdot  q_{j}} 
 =:  S_{1, k}(z)+S_{2, k}(z)+S_{3, k} (z). $$
 We want to prove  the inequality 
\begin{equation} 
\label{f12}\min_{|z| = R_k}|S_{2, k} (z)| > 
\max_{|z| = R_k}( |S_{1, k}(z)|+|S_{3, k} (z)|).
\end{equation}  
As in the previous case, this inequality follows from 
\begin{equation} 
\label{f13}
1-\sum_{j=1}^{k}b_{2n+2}^\frac{-j^2}{2} - \sum_{j=1}^{2n +1 -k}
b_{2n+2}^\frac{-j^2}{2}> 0.
\end{equation}
We have
\begin{equation} 
\label{f14}
1-\sum_{j=1}^{k}b_{2n+2}^\frac{-j^2}{2} - \sum_{j=1}^{2n +1 -k}
b_{2n+2}^\frac{-j^2}{2}> 1 - 2\sum_{j=1}^{n+1}b_{2n+2}^\frac{-j^2}{2} = 0
\end{equation}
by the definition (\ref{f4}) of the constant $b_{2n+2}.$
Using  Rouch\'{e}'s theorem, we obtain that  all the zeros of $P_{2n+1}$ are simple. 
Moreover, the moduli of all zeros of $P_{2n+1}$  are pairwise different.

It remains to consider the case of entire functions.
 Let $f(z)=1 + z+\frac{z^2}{q_{2}}+\frac{z^3}{q_{2}^2q_{3}}
 +\frac{z^4}{q_{2}^3q_{3}^2 q_4} +\ldots$  be an 
entire function.  Suppose that the inequalities  
$|q_k(f)| \geq b_{\infty}$ hold for all $k \geq 2.$ For all $k \in \mathbb{N}$ we 
have
$$f (z)=\sum_{j=0}^{k-1}\frac{z^j}{q_{2}^{j-1}q_{3}^{j-2}\cdot \ldots \cdot  q_{j}}+
\frac{z^k}{q_{2}^{k-1}q_{3}^{k-2} \cdot \ldots \cdot  q_{k}}+\sum_{j=k+1}^{\infty}
\frac{z^j}{q_{2}^{j-1}q_{3}^{j-2}\cdot \ldots \cdot  q_{j}} $$
$$ =:  S_{1, k}(z)+S_{2, k}(z)+S_{3, k} (z). $$
We want to obtain  the inequality 
\begin{equation} 
\label{f12}\min_{|z| = R_k}|S_{2, k} (z)| > 
\max_{|z| = R_k}( |S_{1, k}(z)|+|S_{3, k} (z)|).
\end{equation}  
This inequality follows from 
\begin{equation} 
\label{f13}
1-\sum_{j=1}^{k}b_{\infty}^\frac{-j^2}{2} - \sum_{j=1}^{\infty}
b_{\infty}^\frac{-j^2}{2}> 0.
\end{equation}
We have
\begin{equation} 
\label{f14}
1-\sum_{j=1}^{k}b_{\infty}^\frac{-j^2}{2} - \sum_{j=1}^{\infty}
b_{\infty}^\frac{-j^2}{2}> 1 - 2\sum_{j=1}^{\infty}b_{\infty}^\frac{-j^2}{2} = 0
\end{equation}
by the definition (\ref{f3}) of the constant $b_{\infty}.$
Using  Rouch\'{e}'s theorem, we obtain that  all the zeros of $f$ are simple,
moreover, the moduli of 
all zeros of $f$  are pairwise different.

Theorem \ref{th:mthm1} is proved.

\section{Proof of  Theorem  \ref{th:mthm2}}

At first we prove Theorem  \ref{th:mthm2} (i).   For every $n\in \mathbb{N}$ 
and $c >0$  we consider the following polynomial
\begin{equation} 
\label{f15}
P_{2n, c}(z) = \sum_{k=0}^{2n}c^{k(2n-k)/2} z^k-2c^{n^2/2} z^n.
\end{equation}
Note that, for all $k=0, 1, \ldots, 2n$ the modulus of the $k$-th coefficient
of $P_{2n, c}$ is equal to $c^{k(2n-k)/2} ,$ so that
\begin{equation} 
\label{f15a} |q_k(P_{2n, c})| = \frac{c^{(k-1)(2n-k+1)}}
{c^{(k-2)(2n-k+2)/2}\cdot c^{k(2n-k)/2}} =c 
\end{equation}
for $k=2, 3, \ldots, 2n.$

We observe that $P_{2n, c}(z)=z^{2n}\cdot P_{2n, c}(\frac{1}{z})$, so $P_{2n, c}$ is a 
self-reciprocal polynomial.   We have
$$P_{2n, c}^\prime (z) =2n  z^{2n-1}\cdot P_{2n, c}(\frac{1}{z})-
z^{2n} P_{2n, c}^\prime(\frac{1}{z})\frac{1}{z^2}.$$
Thus, if $P_{2n, c}(1)=0,$ we get
$$P_{2n, c}^\prime(1)= - P_{2n, c}^\prime(1).$$
It means that if $P_{2n, c}(1)=0$ then  $P_{2n, c}^\prime (1)=0,$  so $1$  is a multiple 
root for this polynomial.
Now we consider the equation 
$$P_{2n, c}(1)=\sum_{k=0}^{2n}c^{k(2n-k)/2}-2c^{n^2/2}=0.$$
We rewrite it in the form $2\sum_{k=0}^{n-1}c^{k(2n-k)/2}-c^{n^2/2}=0.$ 
After dividing by $- c^{n^2/2}$ we get
$$1-2\sum_{k=0}^{n-1}c^{-(n-k)^2/2}=0,$$
or, changing the index in the sum: $n-k= j,$ 
$$1-2\sum_{j=1}^{n}c^{-j^2/2}=0. $$
The unique positive root of this equation is  $b_{2n}$ (see (\ref{f4})), 
so we get that the polynomial  $P_{2n, b_{2n}}$ has a multiple root 
and $ |q_k(P_{2n, b_{2n}})|  = b_{2n} $
for all $k=2, 3, \ldots, 2n.$

Let us fix an arbitrary $n\in \mathbb{N}$ and $\varepsilon > 0.$ To prove 
Theorem  \ref{th:mthm2} (ii) we use the polynomial $P_{2n, b_{2n}}.$
Let $Q_{2n+1, d}(z) = P_{2n, b_{2n}}  (1+\frac{z}{d}),$  where $d>0.$
We have $\deg Q_{2n+1, d} = 2n+1,$ and $Q_{2n+1, d}$ has a multiple
root at the point $1$ since $P_{2n, b_{2n}}$ has a multiple
root at that point. Let $ P_{2n, b_{2n}}(z) = \sum_{k=0}^{2n} a_k z^k, $
then
$$Q_{2n+1, d}(z)  =a_{0}+(\frac{a_{0}}{d}+a_{1})z+(\frac{a_{1}}{d}+a_{2})z^2+
\dots +(\frac{a_{2n-1}}{d}+a_{2n})z^{2n}+\frac{a_{2n}}{d}z^{2n+1}.$$
Thus, for all $k = 3, 4, \ldots , 2n$ we have
$$q_{k} (Q_{2n+1, d})=\frac{ (\frac{a_{k-2}}{d}+a_{k-1})^2 }{(\frac{a_{k-1}}{d}+a_{k}) 
(\frac{a_{k-3}}{d}+a_{k-2})}\to q_{k}(P_{2n, b_{2n}}), \  d\to \infty.$$ 
For $k=2$ and $k=2n+1$ we have 
$$q_{2}(Q_{2n+1, d})=\frac{(\frac{a_{0}}{d}+a_{1})^2 }{a_{0}
(\frac{a_{1}}{d}+a_{2}) }\to  q_{2}(P_{2n, b_{2n}}),\  d\to \infty$$
and
$$q_{2n+1} (Q_{2n+1, d})=\frac{(\frac{a_{2n-1}}{d}+a_{2n})^2 }{\frac{a_{2n}}{d}
(\frac{a_{2n-2}}{d}+a_{2n-1}) } \to  \infty,\  d\to \infty.$$
So, for $d$ being large enough we obtain  $|q_k(Q_{2n+1, d})| > b_{2n} -
\varepsilon $  for all $k = 2, 3, \ldots , 2n+1.$

It remains  to prove Theorem  \ref{th:mthm2} (iii). Let us fix an 
arbitrary  $\varepsilon > 0.$ Since $\lim_{n \to \infty} b_{2n}
= b_\infty,$ there exists $n_0 \in \mathbb{N}$ such that
$b_{2n_0} > b_\infty - \varepsilon/3.$ We consider an entire
function of the form
$$f_\varepsilon(z) =  P_{2n_0, b_{2n_0}}(z)  \prod_{j=1}^\infty \left(
1 + \frac{z}{d_j}  \right),$$
where the polynomial $P_{2n_0, b_{2n_0}} $ is defined by (\ref{f15}),
and positive constants $(d_j)_{j=1}^\infty,$ such that $\sum_{j=1}^\infty
\frac{1}{d_j} < \infty,$ will be chosen inductively. We see that $f_\varepsilon$
has a multiple zero at the point $1.$

We know that $|q_k(P_{2n_0, b_{2n_0}})| > b_\infty - \varepsilon/3$ for all $k=2, 3, \ldots, 
2n_0.$ As we have proved above, for  $d_1 >0 $  being large enough for the
polynomial $T_1(z) := P_{2n_0, b_{2n_0}}(z)  (1 + \frac{z}{d_1})$
we have $|q_k(T_1)| >  b_{\infty} - \varepsilon/3 - \varepsilon/4$ for all $k=2, 3, \ldots, 
2n_0 +1.$ We additionally suppose that $d_1 > 2.$ We fix such $d_1,$
and now choose $d_2.$ For all $d_2 >0 $  being large enough for the
polynomial $T_2(z) := T_1(z)  (1 + \frac{z}{d_2})$
we have $|q_k(T_2)| >  b_{\infty} - \varepsilon/3 - \varepsilon/4 - \varepsilon/8$ 
for all $k=2, 3, \ldots, 2n_0 +2.$ We additionally suppose that $d_2 > 2^2.$ Reasoning 
analogously, we construct a sequence of positive constants $(d_j)_{j=1}^\infty,$ 
such that for every $j$ the polynomial $T_j(z) := P_{2n_0, b_{2n_0}}(z)  
\prod_{l=1}^j \left(1 + \frac{z}{d_l}  \right) $ has the property 
$|q_k(T_j)| >  b_{\infty} - \varepsilon/3 - \varepsilon/4 - \varepsilon/8
- \ldots - \varepsilon/2^{j+1} $ for all $k=2, 3, \ldots, 2n_0 +j,$ and
$d_j > 2^j.$ Thus, $\sum_{j=1}^\infty
\frac{1}{d_j} < \infty,$ and $f_\varepsilon$ is an entire function.
We also observe that for every natural $k\geq 2$
$$|q_k(f_\varepsilon)| \geq b_{\infty} - \varepsilon/3 -\sum_{j=1}^\infty 
\frac{\varepsilon}{2^{j+1}} = b_{\infty} - \varepsilon/3 - \varepsilon/2 
> b_{\infty} - \varepsilon.  $$

Theorem \ref{th:mthm2} is proved.

\section{Proof of  Theorem  \ref{th:mthm3}}

We consider a complex polynomial of degree $3$ with nonzero coefficients

$$P_{3, a, b}(z)= 1 + z +\frac{z^2}{a} + \frac{z^3}{a^2 b},\   
a, b \in \mathbb{C}\setminus \{0\},$$
so that $q_2(P_{3, a, b}) =a, q_3 (P_{3, a, b}) =b.$
The polynomial $P_{3, a, b}$ has multiple roots if and only if its
discriminant is equal to zero. We recall that, if
$Q(z) =\alpha z^3 +\beta z^2 + \gamma z +\delta \in 
\mathbb{C}[z],$  $\alpha \neq 0,$ then its discriminant
is equal to $D(Q) =-4 \beta^3 \delta + \beta^2 \gamma^2 
-4 \alpha \gamma^3 +18 \alpha  \beta \gamma \delta 
-27 \alpha^2 \delta^2.$
So we have
$$ D(P_{3, a, b}) =  -\frac{4}{a^3} +\frac{1}{a^2}-\frac{4}{a^2b}
+ \frac{18}{a^3 b} - \frac{27}{a^4 b^2}   = $$
$$-\frac{1}{a^4 b^2}  \left(4ab^2 - a^2b^2+ 4a^2b-18 ab+27\right).$$   
Thus, $P_{3, a, b}$ has multiple roots if and only if 
\begin{equation} 
\label{f16}
4ab^2 - a^2b^2+ 4a^2b-18 ab+27=0.
\end{equation}
Denote by $S := \{(a, b) : a, b \in \mathbb{C} \setminus \{0\},  
4ab^2 - a^2b^2+ 4a^2b-18 ab+27=0 \},$ and
$$c := \sup_{(a, b) \in S}(\min(|a|, |b|).   $$
By Theorem  \ref{th:mthm1} (ii) we have  $c \leq b_4.$
By the definition of $c,$ if $|a| > c$ and $|b| > c$
then all the zeros of  $P_{3, a, b}$ are simple.
We want to prove that $c= \sqrt{9+6\sqrt{3}}.$

We rewrite (\ref{f16}) in the form
\begin{equation} 
\label{f17}
(4 a -a^2) b^2 + (4a^2 -18a) b +27 =0. 
\end{equation}
By our assumption $a \neq 0,$ consider now the case
$a=4.$ Then we have $b = \frac{27}{8}$ and 
$\min(|a|, |b|) = \frac{27}{8} < \sqrt{9+6\sqrt{3}}.$

Let $(a_0, b_0)\in S,$ $a_0 \neq 4,$ and $|a_0| \neq |b_0|.$ Without
loss of generality we suppose that $|a_0| < |b_0|$ since 
(\ref{f16}) is symmetric with respect to $a, b.$
Let $a_0 = r e^{i \alpha},  r > 0, \alpha \in \mathbb{R}.$
For  $\varepsilon >0$ being small enough we denote
by $a_\varepsilon = (r +\varepsilon) e^{i \alpha},$
such that  $a_\varepsilon \neq  4.$ We have $|a_\varepsilon| > |a_0|.$
Then we paste $a_\varepsilon$ in the equation  (\ref{f17})
and find the solution $ b_\varepsilon,$ such that 
$(a_\varepsilon, b_\varepsilon) \in S.$ By the continuity 
reasoning $\lim_{\varepsilon \to 0} b_\varepsilon = b_0,$ so 
for  $\varepsilon >0$ being very small we have
$|b_\varepsilon| > |a_0|.$ Thus, $\min(|a_0|, |b_0|) =
|a_0| < \min(|a_\varepsilon|, |b_\varepsilon|).$ Thus,
we conclude that
$$  c = \sup_{(a, b) \in S,\ |a| = |b|}(\min(|a|, |b|).  $$
Now let $a \in \mathbb{C} \setminus \{0\}, b = a e^{i \gamma}, 
\gamma \in \mathbb{R}.$ We substitute $a, b$ into (\ref{f16})
and get
$$4a^3 e^{2i \gamma}  - a^4 e^{2i \gamma}+ 4a^3 e^{i \gamma}
-18 a^2 e^{i \gamma}+27=0.$$
We are searching for the root of the last equation with the maximal possible
modulus. We rewrite the equation in the form
$$a^4 e^{2i \gamma}  - 4a^3(e^{2i \gamma} +e^{i \gamma}) 
+ 18 a^2 e^{i \gamma} -27 =0, $$
or
$$a^4 e^{2i \gamma}  - 8 a^3 e^{\frac{3i \gamma}{2}}\cos \frac{\gamma}{2} 
+ 18 a^2 e^{i \gamma} -27 =0. $$
We denote by $x$ such complex number that $a = x e^{\frac{-i \gamma}{2}},$
and note that $|a| = |x|,$ we also denote by $\lambda = \cos \frac{\gamma}{2}, $
$\lambda \in [0, 1].$ After substituting in the last equation we get
\begin{equation}
\label{eq3}
 x^4 - 8 \lambda x^3 +18 x^2 -27 =0, 
\end{equation}
and we are searching for the root of the last equation with the maximal possible
modulus. 

At first let us estimate the maximal real positive root. We have
for $\lambda \in [0, 1]$
$$ x^4 - 8 \lambda x^3 +18 x^2 -27 \geq 
 x^4 - 8  x^3 +18 x^2 -27 =(x+1)(x-3)^3,  $$
 so for all $\lambda \in [0, 1]$ the maximal positive root of 
the  equation (\ref{eq3}) is less than or equal to $3.$
 
Now we estimate the minimal real negative root. We have
for $x = -y, y>0,$ and $\lambda \in [0, 1]$
$$ y^4 + 8 \lambda y^3 +18 y^2 -27 \geq 
 y^4  +18 y^2 -27 $$
 $$=(y^2 + 6\sqrt{3} -9) (y+\sqrt{9 + 6\sqrt{3}})
(y-\sqrt{9 + 6\sqrt{3}}),  $$
so for all $\lambda \in [0, 1]$ the maximal modulus of the 
negative root of the equation (\ref{eq3}) is less than or equal 
to $\sqrt{9 + 6\sqrt{3}}.$

Now we consider non-real roots of the equation (\ref{eq3}).
We will use the classical Ferrari method to solve (\ref{eq3}).
For $w \in \mathbb{C}$ we rewrite the equation (\ref{eq3})
in the form
\begin{equation} 
\label{fer}
(x^2 -4\lambda x + w)^2 - \left ((16 \lambda^2 + 2w -18) x^2 -
8 \lambda w x + (w^2 +27) \right)=0.  
\end{equation}
We want to find such $w \in \mathbb{C},$ that the discriminant
of the quadratic expression in the brackets will be zero.
We have the equation
$$ \frac{D}{4} = 16 \lambda^2 w^2 - (w^2 +27) 
(16 \lambda^2 + 2w -18)= $$ 
$$ 16 \lambda^2 w^2 - 16 \lambda^2 w^2 -2 w^3 +
18 w^2 - 27\cdot 16 \lambda^2 -54 w + 27\cdot 18 =0,  $$
or
$$ w^3 -9 w^2 + 27 w -27 -27(8 - 8 \lambda^2) =0.$$
We get the equation
$$ (w-3)^3 = 27\cdot 8(1- \lambda^2) \geq 0. $$
Let us put $w = 3+ 6 \sqrt[3]{1 - \lambda^2}$
and denote by $t =  \sqrt[3]{1 - \lambda^2}, t\in [0, 1],$
so that $ \lambda^2 = 1- t^3, w = 3 +6 t.$ For the second quadratic
expression from (\ref{fer}) we have
$$ (16 \lambda^2 + 2w -18) x^2 -
 \lambda w x + (w^2 +27) = (16 -16 t^3 + 6 +12t -18) x^2 $$
$$- 8 \sqrt{1 - t^3}(3+6t)  x +(9 +36 t +36 t^2 +27) =
(4 -16 t^3 +12t) x^2 -  $$
$$8 \sqrt{1 - t^3}(3+6t)  x +(36 +36 t + 36 t^2)=$$
$$ \frac{4}{1-t} \cdot \Big((1-t)^2(2t+1)^2 x^2 - 6 (1-t)(1+2t)\sqrt{1 - t^3} x
+9(1-t^3)  \Big) $$
$$ = \frac{4}{1-t} \cdot \Big((1-t)(2t+1) x - 3 \sqrt{1 - t^3}  \Big)^2.  $$
We paste this in (\ref{fer})  and get
$$(x^2 -4 \sqrt{1 - t^3} x + 3 + 6t)^2 -\Big(\frac{2}{\sqrt{1-t}}  \Big)^2
\cdot \Big( (1-t)(2t+1) x - 3 \sqrt{1 - t^3} \Big)^2  =0. $$
We write the last equation in the form
$$ \left(x^2 -4 \sqrt{1 - t^3} x + 3 + 6t - 2\sqrt{1-t}(2t+1)x +
6\sqrt{1+t+t^2}\right) \cdot $$
$$ \left(x^2 -4 \sqrt{1 - t^3} x + 3 + 6t + 2\sqrt{1-t}(2t+1)x -
6\sqrt{1+t+t^2}\right) =0,$$
or
$$\left(x^2   - (4 \sqrt{1 - t^3}-2\sqrt{1-t}(2t+1)) x + 
(3 + 6t +6\sqrt{1+t+t^2})\right) \cdot$$
$$\left(x^2   - (4 \sqrt{1 - t^3}+2\sqrt{1-t}(2t+1)) x + 
(3 + 6t -6\sqrt{1+t+t^2})\right) =0.$$
So, all $4$ roots of the equation (\ref{eq3}) are the roots of two 
quadratic equations above, where $t =  \sqrt[3]{1 - \lambda^2}, 
t\in [0, 1].$ Both of these quadratic equations have real coefficients. 
Since $3 + 6t -6\sqrt{1+t+t^2} <0,$ the roots of the second equation 
are real, and we have proved above that the moduli of the roots are less than 
or equal to $\sqrt{9 + 6\sqrt{3}}.$   If the roots of the first equation are real, 
then we have proved that the moduli of the roots are less than or equal to 
$\sqrt{9 + 6\sqrt{3}}.$ If the roots of the first equation are complex and conjugate, 
then their  moduli are equal to $\sqrt{3 + 6t +6\sqrt{1+t+t^2}}.$  It remains to 
check that for all $t\in [0, 1]$ we have
$$\sqrt{3 + 6t +6\sqrt{1+t+t^2}} \leq \sqrt{9 + 6\sqrt{3}},$$
that is obviously valid. So we have proved that
$$c = \sup_{(a, b) \in S}(\min(|a|, |b|) = \sqrt{9 + 6\sqrt{3}}.   $$

Theorem \ref{th:mthm3} (i) is proved.

To prove Theorem \ref{th:mthm3} (ii) we consider the polynomial
$$Q_{3}(z) = 1+ z + \frac{z^2}{\sqrt{9 + 6\sqrt{3}}} -
 \frac{z^3}{(\sqrt{9 + 6\sqrt{3}})^3},$$ 
so that $\deg Q_3 =3,$ $q_2(Q_3)  =  \sqrt{9 + 6\sqrt{3}},$  
$q_3(Q_3)  =  - \sqrt{9 + 6\sqrt{3}},$ whence $|q_2(Q_3)| = |q_3(Q_3)| =  
\sqrt{9 + 6\sqrt{3}}.$ We recall that a complex polynomial $P_{3, a, b} =  
1 + z +\frac{z^2}{a} + \frac{z^3}{a^2 b}$ has multiple roots if and only if 
$4ab^2 - a^2b^2+ 4a^2b-18 ab+27=0$ (see  (\ref{f16})). For the polynomial
$Q_3$ we have $b= -a$ and the condition for having multiple roots takes the form
$$4a^3 - a^4-  4a^3+18 a^2+27=0 \Leftrightarrow 
a^4 -18 a^2 -27 =0.$$
It is easy to check that $a = \sqrt{9 + 6\sqrt{3}} $ is a root of this equation,
so for such $a$ the polynomial $Q_3$ has multiple roots. These multiple
roots can be found explicitly, but the expression is rather cumbersome.

Theorem \ref{th:mthm3} is proved.

\section{Proof of  Theorem  \ref{th:mthm4}}

Let $n \in \mathbb{N} $ be a given integer, and $P_{2n}(z) = 
\sum_{k=0}^{2n} a_k z^k $, $a_k 
\in \mathbb{R} \setminus \{0\}$ for all $k,$  be a real polynomial. 
Suppose that the inequalities  $|q_k(P_{2n})| \geq b_{2n}$ hold for 
all $k = 2, 3, \ldots , 2n.$ For an arbitrary  $\lambda, 0 < \lambda <1,$  
we consider a real polynomial $P_{2n, \lambda}(z) = 
\sum_{k=0}^{2n} a_k \cdot \lambda^{k^2} z^k. $ We have
for $k =2, 3, \ldots, 2n$
$$ q_k(P_{2n, \lambda}) = 
\frac{a_{k-1}^2\cdot \lambda^{2(k-1)^2} }{
 a_{k-2}\cdot \lambda^{(k-2)^2}\cdot 
 a_{k} \cdot \lambda^{k^2}}= \frac{ q_k (P_{2n})}{ \lambda^2}
 > b_{2n}. $$
By Theorem \ref{th:mthm1} (i) the moduli of all zeros of $P_{2n, \lambda}$  
are pairwise different.
So, $P_{2n, \lambda}$ can not have complex conjugate zeros,
whence all the zeros of $P_{2n, \lambda}$ are real. Since
$\lim_{\lambda \to 1} P_{2n, \lambda}(z) = P_{2n}(z),$
and this limit is uniform on the compacts in $\mathbb{C},$ using
the Hurwitz's theorem  we obtain that all the zeros of $P_{2n}$ are
real.

Theorem  \ref{th:mthm4} (i) is proved. Using  analogous reasoning we prove 
Theorem \ref{th:mthm4} (iii)  and Theorem \ref{th:mthm4} (v).

Statements  (ii), (iv) and (vi) in Theorem  \ref{th:mthm4} can be proved
using an analogous reasoning. To prove Theorem  \ref{th:mthm4} (ii), for 
example, we fix an arbitrary $n \in \mathbb{N} $ and  $\varepsilon >0, $
and consider a polynomial $P_{2n, b_{2n}}(z) = \sum_{k=0}^{2n}b_{2n}^{k(2n-k)/2} 
z^k-2b_{2n}^{n^2/2} z^n$ (see (\ref{f15}). We recall that
$|q_k(P_{2n, b_{2n}})| = b_{2n} $ for all $k=2, 3, \ldots, 2n$
(see (\ref{f15a})), and that  $P_{2n, b_{2n}}$ has a double zero at the point
$1.$ Since the sequence of coefficients of $P_{2n, b_{2n}}$ has two
sign changes (all coefficients, except the $n$-th, are positive, and the $n$-th
coefficient is negative), we conclude, using  Descartes' rule of signs,
that $P_{2n, b_{2n}}$ has not more than two positive zeros. So the 
polynomial  $P_{2n, b_{2n}}$  has exactly two positive zeros counting 
multiplicities. Whence,
$$ P_{2n, b_{2n}} \geq 0 \quad \mbox{for all} \quad  x \geq 0.  $$
Now for a small $\delta >0$ we consider a polynomial 
$Q_{2n, \delta}(x) = P_{2n, b_{2n}}(x) + \delta x^n.$
Then we get $Q_{2n, \delta}(x) >0$ for all $x \geq 0.$
We observe that $q_k(Q_{2n, \delta}) = q_k(P_{2n, b_{2n}})$
for all $k=2, 3, \ldots , n-1,  n+3, n+4, \ldots, 2n,$ so that
$|q_k(Q_{2n, \delta})| = b_{2n}$ for all $k=2, 3, \ldots , n-1,  
n+3, n+4, \ldots, 2n.$ We also see that $|q_n(Q_{2n, \delta})| >
|q_n(P_{2n, b_{2n}})| = b_{2n}$ and $|q_{n+2}(Q_{2n, \delta})| >
|q_{n+2}(P_{2n, b_{2n}})| = b_{2n}.$ Since 
$\lim_{\delta \to 0}|q_{n+1}(Q_{2n, \delta})| = 
|q_{n+1}(P_{2n, b_{2n}})| = b_{2n},$ we obtain that 
$|q_{n+1}(Q_{2n, \delta})| > b_{2n} -
\varepsilon $ for $\delta$ being small enough. Thus, 
$|q_k(Q_{2n, \delta})| > b_{2n} -
\varepsilon $  for all $k = 2, 3, \ldots , 2n$
for $\delta$ being small enough. It remains to show that 
the polynomial $Q_{2n, \delta}$ has nonreal zeros. Suppose that all
the zeros of $Q_{2n, \delta}$ are real. Then, since $Q_{2n, \delta}(x) 
>0$ for all $x\geq 0,$ we get that all the zeros of $Q_{2n, \delta}$ are
negative. Then all the coefficients of $Q_{2n, \delta}$ have the same signs,
but we know that all coefficients, except the $n$-th, are positive, and the $n$-th
coefficient is negative for $\delta$ being small enough. Thus, we have proved that
$Q_{2n, \delta}$ has nonreal roots.  

Theorem \ref{th:mthm4} is proved.

{\bf Acknowledgment.} The authors are deeply grateful to the 
referee for suggestions that improved the quality of the text.

 The research of the third author was supported  by  the National 
Research Foundation of Ukraine funded by Ukrainian State budget in frames of 
project 2020.02/0096 ``Operators in infinite-dimensional spaces:  the interplay between 
geometry, algebra and topology''.

\end{document}